\documentclass[a4paper,11pt]{article}

\usepackage{amsmath}
\usepackage{amsthm}
\usepackage{amssymb}

\usepackage{stmaryrd}

\newtheorem{theorem}{Theorem}[section]
\newtheorem{lemma}[theorem]{Lemma}

\newcommand{\R}{{\mathbb R}}

\newcommand{\esperance}{{\mathbb E}}
\newcommand{\Oun}{\mathcal{O}(1)}

\font\grfont=cmr10 scaled \magstep2

\newcommand{\gun}{\hbox{\grfont \char49}}

\newcommand{\dbar}{\overline{\Delta}}
\newcommand{\fbar}{\overline{F}}
\newcommand{\pibar}{\overline{\pi}}
\newcommand{\mum}{\mu_{\scriptscriptstyle M}}
\newcommand{\mud}{\mu_{\scriptscriptstyle \Delta}}
\newcommand{\mubar}{\mu_{\scriptscriptstyle \overline\Delta}}

\newcommand{\mbar}{\overline m}

\newcommand{\pf}{\mathcal{L}}
\newcommand{\partibar}{\overline{\mathcal{M}}}
\newcommand{\parti}{{\mathcal{M}}}
\newcommand{\ukl}{{U^{k}_{l}}}

\renewcommand\phi{\varphi}

\def\supess{\textup{ess\ sup}}

\renewenvironment{proof}[1]
{\noindent{\bf Proof.}\hspace{0.1cm} #1} {\hfill$\blacksquare$\bigskip\bigskip}

\begin{document}

\begin{center}
{\Large{\bf Devroye Inequality for a Class of\\
Non-Uniformly Hyperbolic Dynamical Systems}}
\end{center}

\bigskip

\begin{center}
J.-R. Chazottes$^a$, P. Collet$^a$ 
and B. Schmitt$^{b,}$\footnotemark[1]\ 
\end{center}

\begin{center}
$^a$Centre de Physique Th{\'e}orique,\\
Ecole Polytechnique, CNRS UMR 7644\\
F-91128 Palaiseau Cedex, France\\
emails: {\tt jeanrene@cpht.polytechnique.fr}\\
 {\tt collet@cpht.polytechnique.fr}
\end{center}

\begin{center}
$^b$ D{\'e}partement de Math{\'e}matiques \\
Universit{\'e} de Bourgogne\\
Facult{\'e} des Sciences Mirande\\
BP 138, 21004 Dijon Cedex, France\\
email~: {\tt schmittb@u-bourgogne.fr}
\end{center}

\footnotetext[1]{{\bf Acknowledgments}. 
BS acknowledges the kind hospitality of the CPhT at
Ecole Polytechnique. The authors acknowledge the anonymous referee for
a very careful reading of the paper.}

\begin{abstract}
In this paper, we prove an inequality, which we call "Devroye inequality",
for a large class of non-uniformly hyperbolic dynamical systems $(M,f)$.
This class, introduced by L.-S. Young, includes families of piece-wise hyperbolic maps (Lozi-like maps),
scattering billiards (e.g., planar Lorentz gas), unimodal and H{\'e}non-like maps. 
Devroye inequality provides an upper bound for the variance of
observables of the form $K(x,f(x),\ldots,f^{n-1}(x))$, where $K$ is
{\em any} separately H{\"o}lder continuous function of $n$ variables.
In particular, we can deal with observables which are not Birkhoff averages.
We will show in \cite{CCS} some
applications of Devroye inequality to statistical properties of this
class of dynamical systems.

\bigskip


\end{abstract}

\section{Introduction}

This paper deals with variance estimates for a class of non-uniformly
hyperbolic dynamical systems 
This class was introduced by L.-S. Young in an abstract way. It is strictly larger than Axiom A
since it encompasses families of piece-wise hyperbolic maps, like the Lozi maps; 
scattering billiards, like the planar periodic Lorentz gas; certain quadratic and H{\'e}non maps. 
In this setting, she was able to prove existence of Sinai-Ruelle-Bowen
measures, exponential decay of correlations and the central limit theorem for
H{\"o}lder continuous observables.

Very informally speaking, the strategy successfully carried out by L.-S. Young for the above systems
is to construct a new dynamical system over a horseshoe-like subset of the original system by using
``Markovian'' return times so as to obtain a ``tower Markov map''. Then one reduces this Markov
extension to an ``expanding map'' by quotienting out stable manifolds.
On this reduced system, it is possible to define the transfer operator acting on a suitable function
space giving back H{\"o}lder continuous observables in the original dynamical system. 
The crucial ``parameter'' of this construction is the tail of Markovian return times with respect
to Lebesgue measure, see \cite{YC} for an informal description of this construction.
For the above examples, this tail is exponentially small.
In \cite{young} the existence of a spectral gap is proved for the
transfer operator for the quotiented tower map.
From this follows an exponential decay of correlations for H\"older continuous observables in the 
the original system. 

In the present paper, we prove an inequality which we call ``Devroye inequality''. In the context of i.i.d. random
variables assuming values in a finite set, this inequality was first obtained by L. Devroye in
\cite{devroye}.
This inequality provides an upper estimate for the variance of any H{\"o}lder continuous observable
computed along orbit segments of length $n$, in terms of the sum of
the square of its H{\"o}lder constants.
The two crucial features of this inequality are that it is valid for {\em any} $n$ and for
any separately H\"older continuous observable. In particular it applies to observables which
are not necessarily time-averages of observables.
We will show in \cite{CCS} how to apply Devroye inequality to obtain
statistical properties for this class of dynamical systems.

In the setting of piece-wise expanding maps of the interval, a much stronger inequality holds, namely
an exponential inequality \cite{CMS}.
It immediately implies Devroye inequality for Lipschitz observables. Our strategy to prove
this inequality in the present setting share the same global strategy as in \cite{CMS}, 
that is to exploit the spectral properties of the transfer operator,
in particular its spectral gap.
However, many crucial points have to be handled differently.
In particular, some complications obviously arise due to the fact
that we have to succeed in transferring information from the quotiented tower map back to
the original system. In particular, we have to control the approximations we make to
transform original observables into observables in the quotiented tower map.

Two open issues naturally appear after the present work. The first one concerns the validity
of the exponential inequality, proved in \cite{CMS} for expanding maps
of the interval, in the present setting. 
The second one is about dynamical systems with tails of Markovian return times which are
sub-exponential, in particular polynomial, as in \cite{youngbis}. Basic examples of
such systems are maps of the interval with indifferent fixed points. We are not able
at present to prove Devroye inequality in the setting of \cite{youngbis}. For such
systems, there is no spectral gap for the transfer operator and completely different
techniques seem to be needed.

\bigskip

{\em Outline of the paper}. In Sect. \ref{setup}, we present in a short self-contained way the
class of dynamical systems introduced by L.-S. Young. In Sect. \ref{main} we state our main result,
i.e. Devroye inequality for the variance of separately H\"older continuous observables. 
Sect. \ref{preparation} is devoted to a brief description of the tower Markov map and its
quotiented version. In particular, we recall the spectral theory of the transfer operator.
In Sect. \ref{proof} we prove our main result.

\section{A class of non-uniformly hyperbolic systems}\label{setup}

In this section, we recall the essential features of the abstract class of dynamical
systems in \cite{young} to have a reasonably self-contained presentation
and to fix the notations. For the complete set of assumptions and more details,
we refer to \cite{young}.

Let $M$ be a finite-dimensional, regular and compact, Riemann manifold (endowed with a distance
$d(\cdot,\cdot)$) and let $f$ be a $C^{1+\epsilon}$ diffeomorphism
($\epsilon>0$). We denote by $m$ the Lebesgue measure on $M$.

\bigskip

{\bf Hyperbolic product structure}. 
We assume that there is a set $\Lambda\subset M$ with an hyperbolic product structure in the
following sense. For some $n\geq 1$, there exists a continuous family of $d$-dimensional unstable
disks $\Gamma^u=\{\gamma^u\}$ 
and a continuous family of $(\textup{dim}\ M-d)$-dimensional stable
disks $\Gamma^s=\{\gamma^s\}$ with
$$
\Lambda= (\cup \gamma^u)\cap (\cup \gamma^s)\,.
$$
Recall that an unstable disk $\gamma^u$ is defined by the property that for each $x,x'\in\gamma^u$
$$
\limsup_{n\to\infty}\frac{1}{n} \log d(f^{-n}(x),f^{-n}(x'))<0
$$
while a stable disk $\gamma^s$ is defined via the same condition with forward iterations
of $f$ instead of backward ones.

For $x\in \Lambda$, writing $\gamma^u(x)$ for the element of $\Gamma^u$ containing $x$, we
assume that {\em each $\gamma^u$-disk meets each $\gamma^s$-disk in exactly one point}, and that
the intersection is transversal with the angles bounded away from zero.

We assume that the Lebesgue measure $m$ is compatible with the
hyperbolic structure in the sense that
for every $\gamma\in\Gamma^u$ we have $m_\gamma(\{\gamma\cap \Lambda\})>0$, where $m_\gamma$
is the measure induced by $m$ on $\gamma$.

\bigskip

{\bf Markovian return times}. We assume there are finitely many or countably many
pairwise disjoint subsets $\Lambda_1,\Lambda_2,\ldots\subset \Lambda$,
with a hyperbolic product structure and integers $R_i\geq R_0>1$ with the properties that
\begin{enumerate}
\item $\cup_i \Lambda_i=\Lambda$, modulo zero Lebesgue sets in the unstable direction.
The ``return-time map'' $R:\cup_i \Lambda_i \to \mathbb{Z}^+$ is defined by $R|_{\Lambda_i}=R_i$
(with a slight abuse, $R$\ can be viewed as a Lebesgue almost everywhere defined function on 
$\Lambda$).
\item For each $x\in\Lambda_i$, we have $f^{R_i}(\gamma^s(x))\subset \gamma^s(f^{R_i}(x))$
and $f^{R_i}(\gamma^u(x))\supset \gamma^u(f^{R_i}(x))$.
\item For each $n$ there are at most finitely many $i$'s with $R_i=n$.
\end{enumerate}
These return times are used to construct the ``tower map'' which is the Markov
extension of $(\cup_{j=0}^\infty f^j(\Lambda),f)$, see below.

Thoroughly we will assume exponential tail for Markovian return times. This means that we assume
there are $C>0$ and $\theta<1$ so that for some $\gamma\in\Gamma^u$
\begin{equation}\label{EXPRETURN}
m_\gamma (\{ x\in\Lambda\ | \ R\geq n\}) \leq C\ \theta^n.
\end{equation}

Next we recall two assumptions that we shall explicitly use in the sequel.

\bigskip

{\bf Uniform contraction along $\gamma^s$-disks}. There exist $C>0$
and $0<\alpha<1$, such that for all $x\in\Lambda$, for each $x'\in\gamma^s(x)$,
and all $n\in\mathbb{Z}^+$ we have
\begin{equation}\label{alpha}
d(f^n(x),f^n(x'))\leq C \alpha^n\,.
\end{equation}

The notion of {\bf separation time} plays a central role. Let 
$s_0: \Lambda\times\Lambda\to \mathbb{Z}^+ \cup \{\infty\}$ be such that
\begin{enumerate}
\item $s_0(x,x')=s_0(\tilde{x},\tilde{x}')$ whenever $\tilde{x}\in\gamma^s(x)$
and $\tilde{x}'\in\gamma^s(x')$.
\item For each $n\in \mathbb{Z}^+$, the maximum number of orbits starting from $\Lambda$
that are pairwise separated before time $n$ is finite (where we say that $x$ and $x'$
are separated before time $k$ if $s_0(x,x')<k$). This is related to
condition 3 above. 
\item For $x,x'\in\Lambda_i$ we have $s_0(x,x')\geq R_i + s_0(f^{R_i}(x),f^{R_i}(x'))$.
\item For $x\in\Lambda_i$, $x'\in\Lambda_j$ with $i\neq j$ but $R_i =R_j$ we have
$$
s_0(x,x') < R_i -1\,.
$$
\end{enumerate}

{\bf Backward contraction and distortion along $\gamma^u$-disks}. The separation time
$s_0$ on $\Lambda\times\Lambda$ is such that for all $x\in\Lambda$, each $x'\in\gamma^u(x)$ and
all $0\leq k\leq n<s_0(x,x')$
\begin{enumerate}
\item $d(f^n(x),f^n(x'))\leq C \alpha^{s_0(x,x')-n}$;
\item $$\log\prod_{i=k}^n \frac{\textup{det}Df^u(f^i(x))}{\textup{det}Df^u(f^i(x'))}\leq
C \alpha^{s_0(x,x')-n}\,.
$$
\end{enumerate}
We denoted by $f^u$ the restriction of $f$ to the $\gamma^u$-disks.

\bigskip

{\bf Sinai-Ruelle-Bowen measure}. It is proved in \cite{young} that
$f$ admits a Sinai-Ruelle-Bowen measure supported on
$\cup_{j=0}^\infty f^j(\Lambda)$, which we will be denoted by $\mum$
in the sequel.

\section{Devroye inequality}\label{main}

A real-valued function of $n$ variables $K:M^n\to\mathbb{R}$ is
called {\bf separately} $\eta$-{\bf H{\"o}lder continuous}
if the following H{\"o}lder constants $L_j=L_j(K)$, $1\leq j\leq n$, are all finite
\begin{equation}\label{Holdercoef}
L_j:=
\sup_{x_1,x_2,\ldots,x_{j-1},x_j,x_{j+1},\ldots,x_n} \;
\sup_{\tilde{x}_j\neq x_j}
\end{equation}
$$
\frac{
|K(x_1,\ldots,x_{j-1},x_{j},x_{j+1},\ldots,x_n)-K(x_1,\ldots,x_{j-1},\tilde{x}_{j},x_{j+1},\ldots,x_n)|
}{d(x_j,\tilde{x}_j)^\eta}\;\cdot
$$
It is convenient to define $L_j=0$ for $j>n$ and $L_0=0$.

\bigskip

We can now formulate the main theorem of this paper. It provides, for
any $n\geq 1$, an estimate on the variance of observables of the form
$K(x,f(x),\ldots,f^{n-1}(x))$ where $K$ is any separately H\"older
continuous function.

\begin{theorem}[Devroye inequality for the variance]\label{Main}
Let $(M,f,\mum)$ be the dynamical system defined above. Then, for any
$0<\eta\leq 1$, there exists a constant $D=D(\eta)>0$ such that
for any $n\geq 1$, for any separately $\eta$-H{\"o}lder continuous function $K$ of $n$ variables,
we have 
$$
\int \left( K(x,f(x),\ldots,f^{n-1}(x))-\int K(y,f(y),\ldots,f^{n-1}(y)) d\mum(y)\right)^2 d\mum(x) 
$$
\begin{equation}\label{Devroyeinequality}
\leq D\  \sum_{j=1}^n L_j^2\,.
\end{equation} 
\end{theorem}

\bigskip

Examples of dynamical systems that fit the class of dynamical systems
defined above include Axiom A attractors; piecewise hyperbolic maps
(Lozi-like mappings); billiards with convex scatterers (including
planar periodic Lorentz gases); quadratic maps 
and H{\'e}non-type attractors (for parameter sets with positive
Lebesgue measure). We refer the reader
to \cite{young,BY,YC} for details.

\section{Preparatory notions and results}\label{preparation}

To prove Devroye inequality, we need to use the spectral gap for the transfer operator
proved by L.-S. Young. We use almost the same notations as in \cite{young}.

\subsection{The tower map $(F,\Delta)$}

Let $F:\Delta\circlearrowleft$ be the ``tower map'' as in \cite{young}.
More precisely, we have
$$
\Delta:= \{z:=(x,q): x\in \Lambda; q=0,1,\ldots,R(x)-1\}
$$
and
$$
F(z)=F(x,q):=\left\{
\begin{array}{l}
(x,q+1)\;\;\textup{if}\;\;\; q+1 <R(x)\\
(f^R x, 0)\;\;\;\textup{if}\;\;\;\; q+1 =R(x)\,.\\
\end{array}\right.
$$

There is a projection map $\pi:\Delta\to\bigcup_{j=0}^{\infty} f^j(\Lambda)$
such that $f\circ \pi = \pi\circ F$.
There is an $F$-invariant measure $\mud$ related to $\mum$ via the equation
$\mum:=\mud\circ\pi^{-1}$.

\bigskip

{\bf Markov partition for $F$}. We denote by $\mathcal{M}=\{\Delta_{q,j}\}$ 
the {\em Markov partition for $F$} built explicitly in \cite{young}. 
It is worth thinking of $\Delta$ as a disjoint union of sets
$\Delta_q$ consisting of those pairs $(x,q)\in\Delta$ the
second coordinate of which is $q$. We can picture $\Delta$ 
as a tower and refer to $\Delta_q$ as the $q^{{\scriptscriptstyle th}}$ level
of the tower. In particular, $\Delta_q$ is a copy of $\{x\in\Lambda: R(x)>q\}$.

One needs to slightly modify the definition of the separation time $s_0(\cdot,\cdot)$ defined
above, to make it compatible with the Markov partition. Define, as in \cite{young},
for all pairs $z,z'$ belonging to the same $\Delta_{q,j}$, the number
\begin{equation}\label{newseptime}
s(z,z'):= \textup{the largest}\;n\geq 0\;\textup{such that for all}\;i\leq n
\end{equation}
$$
F^i(z)\;\textup{lies in the same element of}\;\mathcal{M}\,\textup{as}\, F^i(z')\,.
$$

Note that restricted to $\Delta_0\times \Delta_0$
\begin{equation}\label{ss_0}
s(\cdot,\cdot)\leq s_0(\cdot,\cdot)\,.
\end{equation}

The following consequence of the above definitions will be used
repeatedly in the sequel.

\begin{lemma}\label{aviron}
There exists a constant $C>0$ such that for any $y,y'\in\Delta$ such
that there exist an integer $q$ and two points
$\tilde{y},\tilde{y}'\in\Delta$
satisfying $s(\tilde{y},\tilde{y}')\geq q$, $F^q(\tilde{y})=y$,
and $F^q(\tilde{y}')=y'$, then
\begin{equation}\label{dist}
d(\pi(y),\pi(y'))\leq C\ \alpha^{\min(q,s(y,y'))}\; .
\end{equation}
\end{lemma}

\begin{proof}
Without loss of generality, we can assume that
$\tilde{y},\tilde{y}'\in\Delta_0$ and $s(y,y')>0$.
Therefore there exists an integer $m$ such that 
$\tilde{y},\tilde{y}'\in \Delta_{0,m}$.
Let $Z:= \gamma^s(\pi(\tilde{y}))\cap
\gamma^u(\pi(\tilde{y}'))$. Notice that by assumption this intersection
is not empty and consists of exactly one point belonging to some
$\Lambda_p$ since this set has a hyperbolic product structure.
Let $z$ be the unique point in $\Delta_{0,m}$ such that $\pi(z)=Z$.
Since $Z$ is on the local stable manifold of $\pi(\tilde{y})$, it follows
from the Markov property that for all $j\geq 0$, $F^j(z)$ and $F^j(\tilde{y})$ belong to the same
atom of $\parti$. This immediately implies that
$$
s_0(Z,\pi(\tilde{y}'))\geq s(z,\tilde{y}') \geq q + s(y,y')\; .
$$
We now apply the "backward contraction along $\gamma^u$-disks" for $Z\in
\gamma^u(\pi(\tilde{y}'))$ and $n=q$. Using also the previous
inequality we obtain
$$
d(f^q(Z),\pi(y')) \leq C\ \alpha^{s(y,y')}\; .
$$
On the other hand, from the "uniform contraction along
$\gamma^s$-disks", we have
$$
d(f^q(Z), \pi(y))\leq C\ \alpha^q\; .
$$
The result follows from the triangle inequality. 
\end{proof}

\subsection{The quotiented tower map $(\fbar,\dbar)$ and the transfer operator}

Let $\fbar:\dbar\circlearrowleft$ be the (non-invertible) expanding map
obtained by quotienting out the $\gamma^s$-leaves from $\Delta$. The projection
will be denoted by $\pibar: \Delta\to\dbar$, and we shall use the notations 
$\{\dbar_{q}\}$, $\{\dbar_{q,j}\}$, etc. with the obvious meanings.
Notice that $\partibar=\{\dbar_{q,j}\}$ is a Markov partition
for $\fbar$.

\bigskip

Let $\mbar$ be the reference measure on $\dbar$ constructed in \cite{young}. On each
$\gamma\in\Gamma^u$, ${\mbar}_{\gamma}$ is absolutely continuous wrt $m_\gamma$.

Before introducing the suitable Banach space on which will act the transfer operator, we
recall the following facts established in \cite{young}:

\bigskip

{\bf Invariant measure for $\fbar$}. 
The map $\fbar:\dbar\circlearrowleft$ has an invariant
probability measure $\mubar$ of the form $d\mubar= \phi\ d\mbar$, where $\phi$ satisfies
\begin{equation}\label{c0}
c_0^{-1}\leq \phi \leq c_0\quad\textup{for some}\;c_0>0 
\end{equation}
and
\begin{equation}\label{regularitedeladensite}
|\phi(z)-\phi(z')|\leq C\ \alpha^{s(z,z')/2}\quad\forall z,z'\in\dbar_{q,j}
\end{equation}
where $\alpha$ is defined at (\ref{alpha})).
This result of course motivates the choice of the function space.

\bigskip

{\bf Regularity of the Jacobian}. In \cite{young}, it is explained how to give
a ``differentiable structure'' so that one can define the Jacobian $J\fbar=|\textup{det}D\fbar|$.
We have the properties
\begin{equation}\label{plouf}
J\fbar\equiv 1\quad\textup{on}\; \dbar\backslash\fbar^{-1}(\dbar_0)
\end{equation}
and
\begin{equation}\label{regularityofthejacobian}
\left|
\frac{J\fbar(y)}{J\fbar(y')}-1\right|
\leq C \alpha^{s(\fbar(y),\fbar(y'))/2}\quad\forall y,y'\in \dbar_{q,j}\cap \fbar^{-1}\dbar_0\,.
\end{equation}

\bigskip

{\bf Function space}.
For any $\sigma$ such that $\sqrt{\alpha}<\sigma<1$, let
$X_\sigma=\{g:\dbar\to\mathbb{R}, \|g\|<\infty\}$ where the norm
$\|\cdot \|$ 
is defined as follows. Writing $g_{q,j}=g|\{\dbar_{q,j}\}$ and
letting $|\cdot|_p$ denote the $L^p$-norm wrt the reference measure $\mbar$ we set
\[
\| g\| := \|g\|_{\infty} + \|g\|_h
\]
where
\[
\|g\|_{\infty}:=\sup_{q,j}\|g_{q,j}\|_\infty \; ,\quad\|g\|_h :=\sup_{q,j}\|g_{q,j}\|_h
\]
and $\|g_{q,j}\|_{\infty}$ and $\|g_{q,j}\|_h$ are defined by
\[
\|g_{q,j}\|_{\infty}:=|g_{q,j}|_{\infty} e^{-q \varepsilon}
\]
where $\varepsilon>0$ will be chosen  adequately small later on
and
\[
\|g_{q,j}\|_h:=\left(
\supess_{y,y'\in {\scriptscriptstyle{\dbar}_{q,j}}}
\frac{|g(y)-g(y')|}{\sigma^{s(y,y')}} 
\right)\ e^{-q \varepsilon}\,.
\]
It is easy to verify that $(X_\sigma,\|\cdot\|)$ is a Banach space (parametrized by $\varepsilon$).

The transfer operator associated with the dynamical system $\fbar:\dbar\circlearrowleft$ and
reference measure $\mbar$ is then defined by
\[
\mathcal{P}g(y)= \sum_{y':\fbar(y')=y} \frac{g(y')}{J\fbar(y')}\,\cdot
\]

The normalized transfer operator associated to $\mathcal{P}$ is defined as
\begin{equation}\label{normalizedPF}
\pf g(y)= \frac{1}{\phi(y)} \sum_{y':\fbar(y')=y} \frac{\phi(y')}{J\fbar(y')}\ g(y')
\end{equation}
which satisfies $\pf 1 =1$.
The following spectral property of $\mathcal{L}$ is easily derived from
\cite{young}.

\begin{lemma}\label{LYI}
For any $\sigma\in (\sqrt{\alpha},1)$, there exist constants $C>0$ and
$\rho\in (0,1)$ such that, for all $g\in X_\sigma$, and for any
integer $n$, we have
\begin{equation}
\left\|\ \mathcal{L}^n g -\int g\ d\mubar \ \right\| \leq C \rho^n  \| g \| \,.
\end{equation}
\end{lemma}

\section{Proof of Devroye inequality}\label{proof}

{\bf Preparatory approximations for observables}.
Let $K:M^n \longrightarrow \R $ be a separately $\eta$ H\" older
function of $n$ variables.

Let us use the short-hand notations
$$
\esperance (K)=\int \ K (x,\ f(x),\ldots ,f^{n-1}(x))\ d\mum (x)
$$
and
$$
{\textup var} (K)=\int \ (K(x,\ f(x),\ldots ,f^{n-1}(x))-\esperance (K))^2 \ d\mum (x)
$$
A standard computation gives~:
$$
{\textup var}(K)=\frac12\int\ [K(x,\ldots,f^{n-1}(x))-K(x',\ldots ,f^{n-1}(x'))]^2 d\mum(x)\ d\mum(x').
$$ 
Since by construction $\mum=\mud \circ \pi ^{-1}$ we also have 
$$
{\textup var}(K)=
$$
\begin{equation}\label{voiture}
\frac{1}{2}\int
\Bigl [{\widetilde K} (y,\ldots,F^{n-1}(y))-{\widetilde K} (y',\ldots
,F^{n-1}(y'))\Bigr]^2 d\mud (y) d\mud (y')
\end{equation}
where
$$
{\widetilde K} (y_1,\ldots ,y_n)=K(\pi (y_1) ,\ldots ,\pi (y_n))\ .
$$

We now introduce a new piece-wise constant observable $V$ on $\Delta^n$. 
We will write $\parti(x)$ for the atom of the partition $\parti$
containing $x$.

For a fixed integer $p_{0}$ large enough,
we define the integer-valued function $\ell : \mathbb{N} \to \mathbb{N}\cup\{0\}$ by  
\begin{equation}\label{lulu}
\ell(k):=\left\{\begin{array}{ll}
k-1 & \textrm{if}\; k\leq p_0 \log (1+k)\\
p_{0}\big[\log (1+k)\big] & \textup{otherwise}\,.
\end{array}
\right.
\end{equation}

We now define the function $V : \Delta ^n \longrightarrow \R $ as follows.
$$
V(y_1 ,\ldots,y_n):=
$$
\begin{equation}\label{defV}
\inf_{x_1,\ldots,x_n}
\{ 
{\widetilde
  K}(F^{\ell(1)}(x_1),F^{\ell(2)}(x_2),\ldots,F^{\ell(n)}(x_n))\ : 
\end{equation}
$$
F^{k}(x_j) \in\parti(y_{j-\ell(j)+k})\ ,\
k=0,1,\ldots,2\ell(j)\ ,\ j=1,\ldots,n
\}\; .
$$
If the above set is empty, then we set $V(y_1,\ldots,y_n)=0$.

One can remark immediately that $V$ factorizes through $\pibar$ in the sense that~:
$$
V(y_1,\ldots ,y_n)=U\ (\pibar(y_1),\ldots ,\pibar (y_n))
$$
where $U : {\dbar }^n\longrightarrow \R$ is defined as

$$
U(z_1,z_2,\ldots ,z_n):=
$$
\begin{equation}\label{defU}
\inf_{x_1,\ldots,x_n}
\{ 
{\widetilde
  K}(F^{\ell(1)}(x_1),F^{\ell(2)}(x_2),\ldots,F^{\ell(n)}(x_n))\ : 
\end{equation}
$$
F^{k}(x_j) \in\parti(\pibar^{-1}(z_{j-\ell(j)+k}))\ ,\
k=0,1,\ldots,2\ell(j)\ ,\ j=1,\ldots,n
\}\; .
$$
If the above set is empty, then we set $U(z_1,\ldots,z_n)=0$.

We have the following lemma which allows to replace the observable
$\tilde{K}$ by the piece-wise constant observable $V$.

\begin{lemma}\label{pate}
There is constant $C>0$ such that for $p_0$ large enough (see
(\ref{lulu})) we have
$$
\sup_{y\in \Delta} \left| \tilde{K}(y,\ldots ,F^{n-1}(y))-V(y,\ldots
  ,F^{n-1}(y))\right| \leq C\ \sum_{k=1}^{n} \frac{L_k}{k}\,\cdot
$$
\end{lemma}

\begin{proof}
Given $y\in\Delta$, let $x_1,\ldots,x_j,\ldots,x_n$ be a sequence such that for any $1\leq j\leq n$
\begin{equation}\label{pelle}
F^{k}(x_j) \in\parti(F^{j-\ell(j)+k-1}(y))\quad\textup{for}\quad
k=0,1,\ldots,2\ell(j)\; .
\end{equation}
 We have the identity
$$
\tilde{K}(y,\ldots ,F^{n-1}(y))-\tilde{K}(F^{\ell(1)}(x_1),F^{\ell(2)}(x_2),\ldots,F^{\ell(n)}(x_n))=
$$
$$
\sum_{p=-1}^{n-2}
\bigg(
\tilde{K}(y,\ldots ,F^{p+1}(y),F^{\ell(p+3)}(x_{p+3}),\ldots, F^{\ell(n)}(x_n))
$$
$$
-\tilde{K}(y,\ldots ,F^{p}(y),F^{\ell(p+2)}(x_{p+2}),\ldots,F^{\ell(n)}(x_n)\bigg)
$$
where terms with indices out of range are absent.
Therefore using (\ref{Holdercoef}) yields
$$
\left|\tilde{K}(y,\ldots
  ,F^{n-1}(y))-\tilde{K}(F^{\ell(1)}(x_1),F^{\ell(2)}(x_2),\ldots,F^{\ell(n)}(x_n)
)\right|\leq
$$
$$
\sum_{q=1}^{n} L_{q}\ (d(\pi(F^{\ell(q)}(x_q)),f^{q-1}(\pi(y))))^{\eta}\,.
$$
Using (\ref{pelle}) and lemma \ref{aviron} we obtain
$$
d(\pi(F^{\ell(q)}(x_q)),f^{q-1}(\pi(y)))\leq C\ \alpha^{\ell(q)}\; . 
$$
The Lemma follows by choosing $p_0$ large enough in the definition of
$\ell(q)$ in (\ref{lulu}).
\end{proof}

We can now give an approximation of the variance of $K$ in terms of
the piece-wise constant observable $U$ defined on the quotiented
tower $\dbar$.

\begin{lemma}\label{pizza}
We have the following approximation
$$
{\textup var} (K)\leq \int
\Bigl [ U(z,\ldots,{\fbar}^{n-1}(z))-U(z',\ldots,{\fbar}^{n-1}(z'))\Bigr ]^2
d\mubar(z)\ d\mubar(z')
$$
$$
+\quad \Oun \sum ^n_{j=1}\ L^2 _j
$$ 
where $U$ is defined in (\ref{defU}).
\end{lemma}
  
\begin{proof}
Using (\ref{defU}) and the fact that  $\mubar=\mud\circ \pibar^{-1}$ 
we have
$$
\int\ \Bigl (
V(y,\ldots,F^{n-1}(y))-V(y',\ldots,F^{n-1}(y'))\Bigr )^2
\ d\mud (y)\ d\mud (y')
$$ 
$$
=\int\ \Bigl (
U(z,\ldots ,{\fbar}^{n-1}(z))-U(z',\ldots, {\fbar}^{n-1}(z'))\Bigr )^2\ d\mubar(z)\ d\mubar (z')\,.
$$ 

To alleviate notations let us set
$$
U_n (z)=U(z,\ldots ,\fbar^{n-1}(z))
$$
$$
V_n (y)=V(y,\ldots ,F^{n-1}(y)
$$
$$
{\widetilde K}_n(y)={\widetilde K} (y,\ldots,F^{n-1}(y))\,.
$$

By (\ref{voiture}) we have

\begin{equation}
\nonumber
{\textup var}(K)=
\end{equation}
$$
\frac{1}{2} \int\ \Bigl (
{\widetilde K} _n (y)-V_n(y)+V_n (y)+{\widetilde K}_n (y')-V_n
(y')+V_n (y')\Bigr )^2  d\mud(y) \ d\mud (y')
$$
$$
\leq \int\  ( V_n (y)-V_n
(y'))^2\ d\mud(y)\ d\mud (y')\,\, +
$$
$$
\ \int\ \Bigl ( ({\widetilde K}_n
(y)-V_n(y))-({\widetilde K}_n (y')-V_n (y '))\Bigr )^2 \ d\mud (y)\ d\mud (y')\,.
$$ 
Therefore
$$
\textup{var}\ (K)\leq  \int\ ( U_n(z)-U_n (z'))^2\ d\mubar (z)\
d\mubar(z')
$$
$$
+\ 4 \int \ (\tilde{K}_n (y)-V_n(y))^2\ d\mud (y)\,.
$$ 
We now use Lemma \ref{pate} to estimate $|\tilde{K}_n (y)-V_n(y)|$ and
Cauchy-Schwarz inequality, i.e.
$$
\sum_{k=1}^{n} \frac{L_k}{k} \leq \Oun \left( \sum_{k=1}^{n} L_k^2\right)^{1/2}\,.
$$
This ends the proof of Lemma \ref{pizza}. 
\end{proof}

\bigskip
  
\noindent{\bf Martingale procedure}.
As suggested by the previous lemma we will give an upper bound to the integral 

$$
\int \ (U_n (z)-U_n(z'))^2 d\mubar(z)\ d\mubar (z')\,.
$$

To do that we will use the spectral properties of the normalized
transfer operator ${\pf}$ associated to ${\fbar}$, which is defined at
(\ref{normalizedPF}).

We now define an extension of ${\pf}$, also denoted by
${\pf}$. It maps a function $\kappa(x_1,\ldots ,x_n)$ of $n$ variables
on ${\dbar }$ to a function of $(n-1)$ variables, and is given by 

$$
{\pf}\kappa (x_1,\ldots, x_{n-1})=
{1\over \varphi
(x_1)}{\sum_{y:{\fbar}(y)=x_1}}\ {\varphi (y)\over J{\overline F}(y)}\ \kappa(y,x_1,\ldots,x_{n-1}).
$$ 
It is immediate to verify that if the function of one variable $v$ is given by

$$
v(x)=\kappa(x,{\fbar}(x),\ldots ,{\fbar}^{n-1}(x))
$$
then ${\pf} v(x)={\pf}\kappa(x,{\fbar}(x),\ldots ,{\fbar}^{n-2}(x)).$
Moreover if $\kappa$ is a function of $n$ variables and $k<n$ we have

$${\pf}^k \kappa (x_1 ,\ldots ,x_{n-k})=
$$
$$
{1\over \varphi
  (x_1)}\sum_{y:{\fbar}^k(y)=x_1}\ {\varphi (y)\over
  J{\fbar} ^k(y)}\ \kappa(y,\fbar(y),\ldots ,\fbar^{k-1}(y),\ x_1,\ldots
,x_{n-k}).$$ 

For $k\geq n$, we can use the same definition noting that a function
of $n$ variables is also a function of $k$ variables not depending on
the last $(n-k)$ variables.

The extended transfer operator inherits the main properties of the
basic one. In particular the probability measure
$\mubar$ is $\fbar$-invariant, i.e.

\begin{equation}\label{bicorne}
\int\ {\pf}\kappa(x,\ldots ,{\fbar}^{n-2}(x))\ d\mubar (x)=
\int\ \kappa(x,\ldots ,\fbar^{n-1 }(x))\ d\mubar(x)\,.
\end{equation}

The following lemma (reminiscent of a martingale-difference argument) will allow us to use 
later on Lemma \ref{LYI}.

\begin{lemma}\label{lili}
The following identity holds for any $p\geq 0$
$$
\int\ \Bigl
(U_n(y)-U_n(y')\Bigr )^2\ d\mubar(y)\
d\mubar (y')=
$$ 
$$
2\sum_{k=0}^{n-2}\ \int\ \Bigl ( {\pf}^k\
U(y,\ldots,{\fbar}^{n-k-1}(y))-\ {\pf}^{k+1}\ U(
{\fbar} (y),\ldots,{\fbar}^{n-k-1}(y))\Bigr )^2\ d\mubar(y) \, +
$$ 
$$
2\sum_{k=0}^{p} \int\ \left(\pf^k S_n (y) - \pf^{k+1} S_n
    (y)\right)^2\ d\mubar(y) \, +
$$
\begin{equation}\label{finfin}
\int\ \left(\pf^{p+1} S_n (y) - \pf^{p+1} S_n (y')\right)^2\
d\mubar(y)\ d\mubar (y')
\end{equation}
where $S_n(y)=\pf^{n-1} U_n (y)$ is a function which depends only 
on one variable.
\end{lemma}

\begin{proof}
We can write
$$
\int \ \Bigl
(U_n(y)-U_n(y')\Bigr )^2\ d\mubar(y)\
d\mubar (y')=
$$ 
$$
\int \ d\mubar (y)\ d\mubar (y')
$$
$$
\Bigl
(U_n(y)-{\pf} U_{n} ({\fbar}(y))+{\pf} U_{n} (\fbar(y))
-U_n (y')+{\pf} U_{n}({\fbar}(y'))-{\pf}U_{n}({\fbar}(y'))\Bigr)^2
$$ 
$$
=
2\int \ \Bigl ( U_n (y)-{\pf} U_{n}({\fbar}(y))\Bigr ) ^2 \ d\mubar (y) \, +
$$ 
$$
\int \ \Bigl (
{\pf} U_{n} ({\fbar}(y))-{\pf}U_{n}({\fbar}(y'))\Bigr)^2\ 
d\mubar(y)\ d\mubar(y')\; -
$$ 
$$
2\ \left( \int (U_n (y)-{\pf} U_{n}({\fbar}(y)))\ d\mubar(y)\right)^2  \, + \,
2 \ \int \ d\mubar (y)\ d\mubar (y')\ \times
$$
$$
\Bigl
(U_n(y)-{\pf} U_{n} ({\fbar}(y))+{\pf} U_{n} (\fbar(y'))
-U_n (y')\Bigr) 
\Bigl( {\pf} U_{n}({\fbar}(y))-{\pf}U_{n}({\fbar}(y'))\Bigr)\;.
$$
The term before last is equal to zero using the $\fbar$-invariance of
$\mubar$ and (\ref{bicorne}).
Similarly the last term vanishes using the $\fbar$-invariance of
$\mubar$ and the identity
$$
\int U_n (y)\ {\pf} U_{n}({\fbar}(y))\ d\mubar(y)=
\int \big({\pf} U_n (y) \big)^2\ d\mubar(y)
$$
which follows at once from (\ref{bicorne}).
Lemma \ref{lili} follows by iterating this inequality.
\end{proof}

\bigskip

We now need to estimate
$$
\pf^{k+1}U\big(\fbar(y),\ldots,\fbar^{n-k-1}(y)\big)-
\pf^{k}U\big(y,\fbar(y),\ldots,\fbar^{n-k-1}(y)\big)\,.
$$

We will use a decomposition of $U$ into a sum of terms.

For $0\le k\le n-1$ and $0\leq l\leq k$, we define the function $\ukl$ on $\dbar^{3}$ by
$$
\ukl(u,s,y)=
$$
\begin{equation}\label{ukl}
\inf_{E_1(u,l)\ \cap\ E_2(s,l,k)\ \cap \ E_3(y,k,n)}
\{\tilde{K}\big(F^{\ell(1)}(x_1),F^{\ell(2)}(x_2),\ldots,F^{\ell(n)}(x_n)\big)\}
\end{equation}
where with the notation $x_1^n:=(x_1,\ldots,x_n)$
$$
E_1(u,l):=\scriptstyle{\left\{x_1^n\ \big|\ F^q(x_{j})\in \parti(\pibar^{-1}(\fbar^{q+j-\ell(j)-1}(u)))\;\
\mathrm{for} \;\ 0\le q\le 2\ell(j)\;\ \mathrm{and}\;\   1\leq j\leq l\right\}}
$$
$$
E_2(s,l,k):=\scriptstyle{\left\{x_1^n\ \big|\ F^q(x_{j})\in \parti(\pibar^{-1}(\fbar^{j-l-\ell(j-l)+q-1}(s)))\;
\mathrm{for} \;0\le q\le 2\ell(j-l)\;\mathrm{and}\; l+1\le j\le k+1\right\}} 
$$
$$
E_3(y,k,n):=
$$
$$
\scriptstyle{\left\{x_1^n\ \big|\ F^q(x_{j})\in 
\parti(\pibar^{-1}(\fbar^{j-k-\ell(j-k-1)+q-2}(y)))\;\
\mathrm{for} \;\ 0\le q\le 2\ell(j-k-1)\;\ \mathrm{and}\;\ k+2\le j\le n\right\}} \; .
$$
It is convenient to set
$$
E_1(u,0)=\dbar^n\; , \; E_2(s,k+1,k)=\dbar^n\; , \;E_3(y,n-1,n)=\dbar^n\, .  
$$
We define for $0\le l\le k$
\begin{equation}\label{defv}
v_{l}^{k}\big(\xi,y\big)=\int
 \ukl(\xi,s,y)d\mubar(s)
\;. 
\end{equation}
Note that $v_{0}^{k}\big(\xi,y\big)$ does not depend on $\xi$.
We have obviously for $k\geq 1$
$$
U\big(\xi,\fbar(\xi),\ldots,\fbar^{n-1}(\xi)\big)=
U\big(\xi,\fbar(\xi),\ldots,\fbar^{n-1}(\xi)\big)-
v^{k}_{k}\big(\xi,\fbar^{k+1}(\xi)\big)
$$
$$
+v^{k}_{0}\big(\fbar^{k+1}(\xi)\big) 
+\sum_{l=0}^{k-1}\bigg(v^{k}_{l+1}\big(\xi,\fbar^{k+1}(\xi)\big)
-v^{k}_{l}\big(\xi,\fbar^{k+1}(\xi)\big)\bigg)\;.
$$

For $k=0$, the same formula holds without the sum.

By an easy computation one gets 
\begin{equation}\label{telescope}
\pf^{k}U\big(y,\fbar(y),\ldots,\fbar^{n-k-1}(y)\big)=
$$
$$
\pf^{k}\left(U-v^{k}_{k-1}\right)\big(y,\fbar(y),\ldots,\fbar^{n-k-1}(y)\big)
+v^{k}_{0}(\fbar(y)) 
$$
$$
+\sum_{l=0}^{k-1}\pf_{1}^{k-l}w^{k}_{l}(y,\fbar(y))
\end{equation}
where $\pf_{1}$ acts only on the first variable, i.e.
$$
\pf_{1}^{k-l}w^{k}_{l}(y,y')= 
\frac{1}{\phi(y)}\sum_{\fbar^{k-l}(z)=y}
\frac{\phi(z)}{J\fbar^{k-l}(z)}\ w^k_l(z,y')
$$
and 
$w^{k}_{l}$ is defined by
$$
w^{k}_{l}(u,y)=\frac{1}{\phi(u)}\sum_{\fbar^{l}(z)=u}
\frac{\phi(z)}{J\fbar^{l}(z)}
\left(v^{k}_{l+1}(z,y)-v^{k}_{l}(z,y)\right)\;.
$$

\bigskip

\noindent{\bf Regularity estimates}. 
We now estimate the various terms. 
We will use several times the following elementary lemma whose proof
is left to the reader.

\begin{lemma}\label{moulinette}
Let $\Omega_1, \Omega_2$ be two sets and $\Psi$ a real-valued function
on $\Omega_1\times \Omega_2$. Let $\Upsilon_1,\Upsilon_1'$ be two subsets of $\Omega_1$.
Then
$$
\left|
\inf_{\omega_1\in \Upsilon_1,\omega_2\in \Omega_2} \Psi(\omega_1,\omega_2)-
\inf_{\omega_1\in \Upsilon_1',\omega_2\in \Omega_2} \Psi(\omega_1,\omega_2)
\right|
\leq
$$
$$
\sup_{\omega_1\in \Upsilon_1, \omega_1'\in \Upsilon_1', \omega_2\in \Omega_2}
\left|
\Psi(\omega_1,\omega_2)-\Psi(\omega_1',\omega_2)
\right|\;.
$$
\end{lemma}
To apply this lemma we will use the following sequence of sets
$$
\mathcal{E}(u,m):=\scriptstyle{\left\{x\in\Delta\ \big|\ F^q(x)\in
\parti(\pibar^{-1}(\fbar^{q+m-\ell(m)-1}(u)))\,
\mathrm{for} \,0\le q\le 2\ell(m)\right\}}
$$
where $u\in\dbar$, $m$ is an integer. It is useful to observe that
$$
E_1(u,l)=\vartimes_{j=1}^l \mathcal{E}(u,j)\ , \
E_2(s,l,k)=\vartimes_{j=l+1}^{k+1} \mathcal{E}(s,j-l)
$$
$$
E_3(y,k,n)= \vartimes_{j=k+2}^{n} \mathcal{E}(y,j-k-1)\,.
$$

We denote by $\textup{diam}(M)$ the diameter of $M$. 

\bigskip

The first term we have to estimate is 
$$
\sup_{\xi\in\dbar}\bigg|U\big(\xi,\fbar(\xi),\ldots,\fbar^{n-1}(\xi)\big)-
v^{k}_{k}\big(\xi,\fbar^{k+1}(\xi)\big)\bigg|
=
$$
$$
\sup_{\xi\in\dbar}\bigg|U\big(\xi,\fbar(\xi),\ldots,\fbar^{n-1}(\xi)\big)-
\int U^{k}_{k}\big(\xi,\fbar^{\ell(k+1)}(s),\fbar^{k+1}(\xi)\big)\ d\mubar(s)\bigg|
$$
where we have used the invariance of the measure. 
We know apply lemma \ref{moulinette} by taking
$$
\Omega_1=\Delta^{n-k}\, , \, \Omega_2=\vartimes_{j=1}^{k} \mathcal{E}(\xi,j)
$$
$$
\Upsilon_1=\vartimes_{p=k+1}^{n} \mathcal{E}(\xi,p)\, ,\,
\Upsilon_1'=\mathcal{E}(s,1)\times\vartimes_{j=k+2}^{n}\ \mathcal{E}(\fbar^{k+1}(\xi),j-k-1)
$$
and
$$
\Psi(\omega_1,\omega_2)=
\tilde{K}\big(F^{\ell(1)}(x_1),F^{\ell(2)}(x_2),\ldots,F^{\ell(n)}(x_n)\big)
$$
where
$$
\omega_1=(x_{k+1},\ldots,x_n)\; ,\; \omega_2=(x_1,\ldots,x_{k})\,.
$$
We have
$$
|\Psi(\omega_1,\omega_2)-\Psi(\omega_1',\omega_2)|\leq
$$
$$
\sum_{p=k}^{n-1}
\Big|
\tilde{K}\left(F^{\ell(1)}(x_1),\ldots,F^{\ell(p)}(x_p), F^{\ell(p+1)}(x_{p+1}'), F^{\ell(p+2)}(x_{p+2}'), \ldots,
F^{\ell(n)}(x_{n}')\right)
$$
$$
-
\tilde{K}\left(F^{\ell(1)}(x_1),\ldots,F^{\ell(p)}(x_p), F^{\ell(p+1)}(x_{p+1}), F^{\ell(p+2)}(x_{p+2}'), \ldots,
F^{\ell(n)}(x_{n}')\right)
\Big|
$$
For $p=k$, we get the upper bound 
$L_{k+1} (\textup{diam}{M})^\eta$
by using \eqref{Holdercoef}.
For $p\geq k+1$, we apply Lemma \ref{aviron} with
$y=F^{\ell(p+1)}(x_{p+1})$, $y'=F^{\ell(p+1)}(x_{p+1}')$, $q=\ell(p-k)$,
$\tilde{y}= F^{\ell(p+1)-\ell(p-k)}(x_{p+1})$, $\tilde{y}'= F^{\ell(p+1)-\ell(p-k)}(x_{p+1}')$,
observing that $s(y,y')\geq \ell(p-k)$ (this follows from the definition of the sets $\mathcal{E}$).
We finally obtain
$$
|\Psi(\omega_1,\omega_2)-\Psi(\omega_1',\omega_2)|\leq
$$
$$
L_{k+1} (\textup{diam}(M))^\eta + C^{\eta} \sum_{p=k+1}^{\infty} L_{p+1} \alpha^{\eta\ell(p-k)}\,.
$$
Hence 
\begin{equation}\label{pollene}
\sup_{\xi\in\dbar}\bigg|U\big(\xi,\fbar(\xi),\ldots,\fbar^{n-1}(\xi)\big)-
v^{k}_{k}\big(\xi,\fbar^{k+1}(\xi)\big)\bigg|
\leq B_{k+1}
\end{equation}
where, for any $q\in\mathbb{N}$ 
\begin{equation}\label{pt}
B_q := L_{q} (\textup{diam}(M))^\eta + C^{\eta} \sum_{p=q}^{\infty}
L_{p+1} \alpha^{\eta\ell(p-q+1)}\,.
\end{equation}

\bigskip

To estimate $v^{k}_{l}-v^{k}_{l-1}$, we use, as in \cite{CMS}, the invariance of the measure $\mubar$ to write
$$
v^{k}_{l}(\xi,y)-v^{k}_{l-1}(\xi,y)=
$$
$$
\int\ \left(\ukl(\xi,\fbar^{\ell(l+1)}(\fbar(s)),y)-
U^{k}_{l-1}(\xi,\fbar^{\ell(l)}(s),y)\right)d\mubar(s)\;.
$$
To estimate the integrand, we apply lemma \ref{moulinette} by taking
$$
\Omega_1=\Delta^{n-k+l-2}\, , \, 
\Omega_2=\vartimes_{j=1}^{l-1} \mathcal{E}(\xi,j)\times
\vartimes_{j=k+2}^{n} \mathcal{E}(y,j-k-1)
$$
$$
\Upsilon_1=\mathcal{E}(\xi,l) \times \vartimes_{j=l+1}^{k+1} \mathcal{E}(\fbar(s),j-l)\, ,\,
\Upsilon_1'=\vartimes_{j=l}^{k+1}\mathcal{E}(s,j-l+1)
$$
and
$$
\Psi(\omega_1,\omega_2)=
\tilde{K}\big(F^{\ell(1)}(x_1),F^{\ell(2)}(x_2),\ldots,F^{\ell(n)}(x_n)\big)
$$
where
$$
\omega_1=(x_{l},\ldots,x_{k+1})\; ,\; \omega_2=(x_1,\ldots,x_{l-1},x_{k+2},\ldots,x_n)\,.
$$
We have
$$
|\Psi(\omega_1,\omega_2)-\Psi(\omega_1',\omega_2)|\leq
$$
$$
\sum_{p=l}^{k+1}
\Big|
\tilde{K}\left(F^{\ell(1)}(x_1),\ldots,F^{\ell(p)}(x_p'),
  F^{\ell(p+1)}(x_{p+1}'), \ldots, \right.
$$
$$
\left. F^{\ell(k+1)}(x_{k+1}'), F^{\ell(k+2)}(x_{k+2}), \ldots,
F^{\ell(n)}(x_{n})\right)
$$
$$
-
\tilde{K}\left(F^{\ell(1)}(x_1),\ldots,F^{\ell(p)}(x_p),
  F^{\ell(p+1)}(x_{p+1}'), \ldots,\right.
$$
$$
\left. F^{\ell(k+1)}(x_{k+1}'), F^{\ell(k+2)}(x_{k+2}), \ldots,
F^{\ell(n)}(x_{n})\right)
\Big|
$$
For $p=l$, we get the upper bound $L_{l} (\textup{diam}{M})^\eta$
by using \eqref{Holdercoef}.
For $p\geq l+1$, we apply Lemma \ref{aviron} with the two points
$F^{\ell(p)}(x_{p})$, $F^{\ell(p)}(x_{p}')$, $q=\ell(p-l)$,
$\tilde{y}= F^{\ell(p)-\ell(p-l)}(x_{p})$, $\tilde{y}'= F^{\ell(p)-\ell(p-l)}(x_{p}')$,
observing that $s(F^{\ell(p)}(x_{p}), F^{\ell(p)}(x_{p}'))\geq \ell(p-l)$.
We finally obtain
$$
|\Psi(\omega_1,\omega_2)-\Psi(\omega_1',\omega_2)|\leq B_l
$$
where $B_l$ is defined at (\ref{pt}).

It follows that for any $1\leq l \leq k$ 
\begin{equation}\label{lmlm1}
\sup_{\xi,y}\bigg|v^{k}_{l}(\xi,y)-v^{k}_{l-1}(\xi,y)\bigg|\le B_l\,.
\end{equation}
This immediately implies for any $0\le l\le k-1$
\begin{equation}\label{udt}
\big|w^{k}_{l}\big|_{\infty}\le B_{l+1}\;.
\end{equation}

We now have to estimate the regularity of $w^{k}_{l}$ with respect to
its first variable. This is the content of the
following lemma.

\begin{lemma}\label{wreg}
There is a constant $C>0$ such that for any $z$, $z'$ and $y$ in
$\dbar$, for any integers $k,l$ with $0\leq l\leq k-1$, and for any
separately $\eta$-H\"older
continuous observable $K$,
$$
\big|w^{k}_{l+1}(z,\fbar(y))-w^{k}_{l}(z',\fbar(y))\big|\le 
$$
$$
C \ \alpha^{\eta s(z,z')/2}\
\left(B_l + (\textup{diam}(M))^{\eta} 
\sum_{j=0}^{l-1} \alpha^{\eta (l-j)/2}\ L_{j}\right)\,.
$$
\end{lemma}

\begin{proof}
It is convenient to distinguish two cases. The first case
corresponds to $s(z,z')=0$. We then use the estimate (\ref{udt}) and the
result follows. We now consider the case $s(z,z')>0$. Using the Markov
property of the map $\fbar$ on $\dbar$, we can write in this case
$$
w^{k}_{l}(z,y)-w^{k}_{l}(z',y)=
\sum_{\zeta\in\bigvee_{j=1}^{l}\fbar^{-j}\partibar}\;\:
\gun_{\zeta\cap \fbar^{-l}(z)}(\xi)\
\gun_{\zeta\cap \fbar^{-l}(z')}(\xi')
$$
\begin{equation}\label{ladif}
\left(\frac{\phi(\xi)}{\phi(z)J\fbar^{l}\!(\xi)}
\left(v^{k}_{l}(\xi,y)-v^{k}_{l-1}(\xi,y)\right)
-\frac{\phi(\xi')}{\phi(z')J\fbar^{l}\!(\xi')}
\left(v^{k}_{l}(\xi',y)-v^{k}_{l-1}(\xi',y)\right)
\right)\,.
\end{equation}

We first observe that using properties (\ref{regularitedeladensite})
and (\ref{regularityofthejacobian}), and the fact that $s(\xi,\xi')\ge s(z,z')$,
we get for some uniform constant $C_{1}>0$
$$
\left|\frac{\phi(\xi)}{\phi(z)J\fbar^{l}(\xi)}
-\frac{\phi(\xi')}{\phi(z')J\fbar^{l}(\xi')}\right|
\le C_{1}\ \alpha^{s(\xi,\xi')/2} \frac{\phi(\xi)}{\phi(z)J\fbar^{l}(\xi)}\;\cdot
$$
Therefore, using the estimate (\ref{lmlm1}) we obtain
\begin{equation}\label{unmorceau}
\left|\left(\frac{\phi(\xi)}{\phi(z)J\fbar^{l}(\xi)}
-\frac{\phi(\xi')}{\phi(z')J\fbar^{l}(\xi')}\right)
\left(v^{k}_{l}(\xi,y)-v^{k}_{l-1}(\xi,y)\right)
\right|
$$
$$
\le (1+C_{1})\ \alpha^{s(z,z')/2}\ B_{l+1} \
\frac{\phi(\xi)}{\phi(z)J\fbar^{l}(\xi)}\;\cdot
\end{equation}
It remains to estimate 
$$
\frac{\phi(\xi')}{\phi(z')J\fbar^{l}(\xi')}\
\bigg[\left(v^{k}_{l}(\xi,y)-v^{k}_{l-1}(\xi,y)\right)
-\left(v^{k}_{l}(\xi',y)-v^{k}_{l-1}(\xi',y)\right)\bigg]
$$
for $\xi$ and $\xi'$ in the same atom of 
$\bigvee_{j=0}^{l+s(z,z')-1}\fbar^{-j}\partibar$. Coming back to the definition
(\ref{defv})  of $v^{k}_{l}$, we get
$$
v_{l}^{k}\big(\xi,y\big)-v_{l}^{k}\big(\xi',y\big)
=\int \left(\ukl(\xi,s,y)-\ukl(\xi',s,y)\right)\ d\mubar(s)\;.
$$
We are going to prove that if $\xi$ and $\xi'$ belong to 
the same atom of $\bigvee_{j=0}^{l-1}\fbar^{-j}\partibar$, 
with $\fbar^l(\xi)=z$ and $\fbar^l(\xi')=z'$, we have
$$
\sup_{s,y}\left|\ukl(\xi,s,y)-\ukl(\xi',s,y)\right|\le
$$
\begin{equation}\label{quiche}
C \ (\textup{diam}(M))^{\eta}\ \alpha^{\eta s(z,z')/2}\
\sum_{j=1}^{l} \alpha^{\eta (l-j-1)/2}\ L_{j}
\end{equation}
where $C>0$ is a uniform constant.

First observe that if $s(z,z')\geq \ell(l)$ then it follows
immediately from definition (\ref{ukl}) that 
$\ukl(\xi,s,y)=\ukl(\xi',s,y)$. Hence the estimate is true in this
case.

Now assume that $s(z,z')< \ell(l)$. Let $p_*=p_*(l)$ be the largest
integer such that $p_* + \ell(p_*) < l$. 

Observe that for any $1\leq j\leq l-\ell(l)$ we have
$\mathcal{E}(\xi,j)=\mathcal{E}(\xi',j)$ since by assumption $\xi$ and $\xi'$ belong to 
the same atom of $\bigvee_{j=0}^{l-1}\fbar^{-j}\partibar$.

We now apply lemma \ref{moulinette} by taking
$$
\Omega_1=\Delta^{l-p_*}\, , \, 
\Omega_2=\vartimes_{j=1}^{p_*} \mathcal{E}(\xi,j)
\times \vartimes_{j=l+1}^{k+1} \mathcal{E}(s,j)
\times\vartimes_{j=k+2}^{n} \mathcal{E}(y,j-k-1)
$$
$$
\Upsilon_1=\vartimes_{j=p_{*}+1}^{l+1} \mathcal{E}(\xi,j)
\, ,\,
\Upsilon_1'=\vartimes_{j=p_{*}+1}^{l+1} \mathcal{E}(\xi',j)
$$
and
$$
\Psi(\omega_1,\omega_2)=
\tilde{K}\big(F^{\ell(1)}(x_1),F^{\ell(2)}(x_2),\ldots,F^{\ell(n)}(x_n)\big)
$$
where
$$
\omega_1=(x_{p_{*}+1},\ldots,x_{l})\; ,\; \omega_2=(x_1,\ldots,x_{p_{*}},x_{l+1},\ldots,x_n)\,.
$$
We have
$$
|\Psi(\omega_1,\omega_2)-\Psi(\omega_1',\omega_2)|\leq
$$
$$
\sum_{p=p_{*}+1}^{l}
\Big|
\tilde{K}\left(F^{\ell(1)}(x_1),\ldots,F^{\ell(p)}(x_p'),
  F^{\ell(p+1)}(x_{p+1}'), \ldots, \right.
$$
$$
\left. F^{\ell(k+1)}(x_{k+1}'), F^{\ell(k+2)}(x_{k+2}), \ldots,
F^{\ell(n)}(x_{n})\right)
$$
$$
-
\tilde{K}\left(F^{\ell(1)}(x_1),\ldots,F^{\ell(p)}(x_p),
  F^{\ell(p+1)}(x_{p+1}'), \ldots,\right.
$$
$$
\left. F^{\ell(k+1)}(x_{k+1}'), F^{\ell(k+2)}(x_{k+2}), \ldots,
F^{\ell(n)}(x_{n})\right)
\Big|
$$
We apply Lemma \ref{aviron} with the two points
$F^{\ell(p)}(x_{p})$,
$F^{\ell(p)}(x_{p}')$,
$q=\ell(p)$,
$\tilde{y}=x_{p}$,
$\tilde{y}'= x_{p}'$,
observing that $s(x_{p}, x_{p}')\geq \min(l+s(z,z'), p+\ell(p)) -p$.
We finally obtain
$$
|\Psi(\omega_1,\omega_2)-\Psi(\omega_1',\omega_2)|\leq 
\sum_{p=p_{*}+1}^{l} L_p \ \alpha^{\eta (\min(\ell(p),l-p+s(z,z')))}\,.
$$
We claim that there exists a number $c_0$ such that for any $l,z,z'$
and $p$ such that $l\geq p \geq p_* +1$, one has 
$\min(\ell(p),l-p+s(z,z'))\geq s(z,z')/2 + (l-p)/2 - c_0$.
This is obvious if $l-p+ s(z,z') \leq \ell(p)$. From the definition
of $p_*$ it follows that there exists a constant $c_1>0$ such that
$\ell(p_*) \geq \ell(l) -c_1$. This implies (since $l\geq p\geq p_*$)
$\ell(p)\geq \ell(p_*)\geq \ell(l)-c_1 \geq s(z,z')/2 + \ell(l)/2-c_1$.
This follows from the monotonicity of $\ell$ and the assumption
$s(z,z') < \ell(l)$. On the other hand, from the definition of $p_*$
we have $\ell(l)\geq \ell(p_* +1) \geq l -p_* -1 \geq l-p$.
Therefore we get the estimate (\ref{quiche}).

It immediately follows from the definition that
\begin{equation}\label{salade}
\bigg|v_{l}^{k}\big(\xi,y\big)-v_{l}^{k}\big(\xi',y\big)\bigg|\le 
C^\eta\ (\textup{diam}(M))^{\eta} \ \alpha^{\eta s(z,z')/2}\
\sum_{j=1}^{l} \alpha^{\eta (l-j-1)/2}\ L_{j}\;.
\end{equation}

\bigskip

Using the estimate (\ref{unmorceau}) we get
$$
\left|
\frac{\phi(\xi)}{\phi(z)J\fbar^{l}(\xi)}
\left(v^{k}_{l}(\xi,y)-v^{k}_{l-1}(\xi,y)\right)
-\frac{\phi(\xi')}{\phi(z')J\fbar^{l}(\xi')}
\left(v^{k}_{l}(\xi',y)-v^{k}_{l-1}(\xi',y)\right)
\right|
$$
$$
\le 
(1+C_{1})\ \alpha^{s(z,z')/2}\ B_{l+1} \
\frac{\phi(\xi)}{\phi(z)J\fbar^{l}(\xi)}
$$
$$
+\; C^\eta\ (\textup{diam}(M))^{\eta} \ \alpha^{\eta s(z,z')/2}\
\frac{\phi(\xi')}{\phi(z')J\fbar^{l}(\xi')}\ 
\sum_{j=1}^{l} \alpha^{\eta (l-j-1)/2}\ L_{j}\;.
$$
The lemma follows from relation (\ref{ladif}) by summing over $\zeta$
and using the identity $\pf 1=1$.
\end{proof}

\bigskip

It follows immediately from the estimate (\ref{udt}) and Lemma
\ref{wreg} that for fixed $y$, as a function of $u$, 
$w^{k}_{l}(u,\fbar(y))$ belongs to the space $X_{\sigma}$, where
$\sigma=\alpha^{\eta/2}$, with an $X_\sigma$-norm satisfying uniformly in $y$ and $k$
$$
\|w^{k}_{l}(\cdot,\fbar(y))\|\le
\Oun \ \left(B_{l+1} + (\textup{diam}(M))^{\eta} 
\sum_{j=1}^{l} \alpha^{\eta (l-j)/2}\ L_{j}\right)\,.
$$

Using lemma \ref{LYI}, we get for some constants $C>0$ and  $0<\rho<1$
independent of $K$, $l$ and $k$,  
\begin{equation}\label{estimw}
\|\pf_{1}^{k-l}w^{k}_{l}(\cdot,\fbar(y))-a_{k,l}(\fbar(y))\|
\le
C \ \Gamma_{l}^{k}
\end{equation}
where 
\begin{equation}\label{gamma}
\Gamma_{l}^{k}=\rho^{k-l} \ 
\left(B_{l+1} + (\textup{diam}(M))^{\eta} 
\sum_{j=1}^{l} \alpha^{\eta (l-j)/2}\ L_{j}\right)
\end{equation}
and
$$
a_{k,l}(y')=\int\ w_l^k(u,y')\ d\mubar(u)\,.
$$

\bigskip

\noindent{\bf Final estimates}.
We start by estimating the first term in (\ref{finfin}).
We observe that
$$
\pf^{k+1}U\big(\fbar(y),\ldots,\fbar^{n-k-1}(y)\big)-
\pf^{k}U\big(y,\fbar(y),\ldots,\fbar^{n-k-1}(y)\big)
$$
$$
=\frac{1}{\phi(\fbar(y))}\sum_{\fbar(u)
=\fbar(y)}\frac{\phi(u)}{J\fbar(u)}
\times
$$
$$
 \bigg[\pf^{k}U\big(u,\fbar(y),\ldots,\fbar^{n-k-1}(y)\big)-
\pf^{k}U\big(y,\fbar(y),\ldots,\fbar^{n-k-1}(y)\big)\bigg]
$$
where we have used the fact that $\pf 1=1$.

We obtain, using equations (\ref{telescope}) and (\ref{pollene}), and
observing that $v_0^k(\fbar(u))=v_0^k(\fbar(y))$, the following
estimate
\begin{equation}\label{borneu}
\int \left(
\pf^{k+1}U\big(\fbar(y),\ldots,\fbar^{n-k-1}(y)\big)-
\pf^{k}U\big(y,\fbar(y),\ldots,\fbar^{n-k-1}(y)\big)\right)^{2}d\mubar(y)
$$
$$
\le \Oun \ B_{k+1}^{2}+
\Oun\int \left(\frac{1}{\phi(\fbar(y))}
\sum_{\fbar(u)=\fbar(y)}\frac{\phi(u)}{J\fbar(u)}\right.
$$
$$
\left.\left(
\sum_{l=0}^{k-1}\big(\pf_{1}^{k-l}w^{k}_{l}(u,\fbar(y))-
\pf_{1}^{k-l}w^{k}_{l}(y,\fbar(y))\big)
\right)
\right)^{2}d\mubar(y)\;.
\end{equation}
Since $\fbar(u)=\fbar(y)$, we have 
$$
 \left(\frac{1}{\phi(\fbar(y))}
\sum_{\fbar(u)=\fbar(y)}\frac{\phi(u)}{J\fbar(u)}
\sum_{l=0}^{k-1}\big(\pf_{1}^{k-l}w^{k}_{l}(u,\fbar(y))-
\pf_{1}^{k-l}w^{k}_{l}(y,\fbar(y))\big)
\right)^{2}
$$
$$
=\left(\frac{1}{\phi(\fbar(y))}
\sum_{\fbar(u)=\fbar(y)}\frac{\phi(u)}{J\fbar(u)}
\sum_{l=0}^{k-1}\big(\pf_{1}^{k-l}w^{k}_{l}(u,\fbar(y))-
a_{k,l}(\fbar(u))\big)\right.
$$
$$
\left.-\frac{1}{\phi(\fbar(y))}
\sum_{\fbar(u)=\fbar(y)}\frac{\phi(u)}{J\fbar(u)}
\sum_{l=0}^{k-1}\big(\pf_{1}^{k-l}w^{k}_{l}(y,\fbar(y))-
a_{k,l}(\fbar(y))\big)\right)^{2}
$$
$$
\le 2  
 \left(\frac{1}{\phi(\fbar(y))}
\sum_{\fbar(u)=\fbar(y)}\frac{\phi(u)}{J\fbar(u)}
\sum_{l=0}^{k-1}\big(\pf_{1}^{k-l}w^{k}_{l}(u,\fbar(y))-
a_{k,l}(\fbar(u))\big)\right)^{2}
$$
\begin{equation}\label{laderniere}
+2\left(
\sum_{l=0}^{k-1}\big(\pf_{1}^{k-l}w^{k}_{l}(y,\fbar(y))-
a_{k,l}(\fbar(y))\big)\right)^{2}\;.
\end{equation}
We now estimate separately the integral of each term. 

We define the integer valued function $q(y)$ by
$$
q(y)=q \;\textrm{if}\;y\in\dbar_{q,j}\;.
$$
We have from (\ref{estimw}) and the definition of the norm in
$X_{\sigma}$ the following estimate uniform in $\fbar(y)$
\begin{equation}\label{youpla}
\bigg|\pf_{1}^{k-l}w^{k}_{l}(u,\fbar(y))-
a_{k,l}(\fbar(y))\bigg|\le \Oun \ e^{\varepsilon q(u)}\ \Gamma_l^k\;.
\end{equation}
We have, since $\fbar(u)=\fbar(y)$, 
$$
\int
\left(\frac{1}{\phi(\fbar(y))}
\sum_{\fbar(u)=\fbar(y)}\frac{\phi(u)}{J\fbar(u)}
\sum_{l=0}^{k-1}\big(\pf_{1}^{k-l}w^{k}_{l}(u,\fbar(y))-
a_{k,l}(\fbar(u))\big)\right)^{2}d\mubar(y)
$$
$$
\le \Oun \left( \sum_{\ell=1}^{k-1} \Gamma^{k}_{l}\right)^{2}
\int 
\left(\frac{1}{\phi(\fbar(y))}
\sum_{\fbar(u)=\fbar(y)}\frac{\phi(u)}{J\fbar(u)}\ e^{\varepsilon q(u)}
\right)^{2}d\mubar(y)
$$
$$
=\Oun \left(\sum_{\ell=1}^{k-1} \Gamma^{k}_{l}\right)^{2}
\int 
\left(\frac{1}{\phi(y)}
\sum_{\fbar(u)=y}\frac{\phi(u)}{J\fbar(u)}\ e^{\varepsilon q(u)}
\right)^{2}d\mubar(y)
$$
by the invariance of the measure $\mubar$.
Using Cauchy-Schwarz inequality and the property $\pf 1=1$, the last
integral is bounded above by
$$
\int 
\frac{1}{\phi(y)}
\sum_{\fbar(u)=y}\frac{\phi(u)}{J\fbar(u)}\ e^{2\varepsilon q(u)}
d\mubar(y)=
\int e^{2\varepsilon q(y)}d\mubar(y)\;.
$$
The integral of the last term in the estimate (\ref{laderniere}) is bounded by the same quantity.
By the invariance of the measure $\mubar$, we have
$$
\int  
e^{2\varepsilon q(y)}d\mubar(y)
=\sum_{j}\sum_{l=0}^{R_{j}-1}e^{2\varepsilon l}
\mubar\big(\dbar_{l,j}\big)
=\sum_{j}\mubar\big(\dbar_{0,j}\big)\sum_{l=0}^{R_{j}-1}e^{2\varepsilon l}
$$
$$
\le \frac{1}{e^{2\varepsilon}-1}\sum_{j}\mubar\big(\dbar_{0,j}\big)\ e^{2\varepsilon R_{j}}\;.
$$
Since $\phi$ is bounded on $\dbar_{0}$, see (\ref{c0}), we get 
$$
\sum_{j}\mubar\big(\dbar_{0,j}\big)
e^{2\varepsilon R_{j}}\le\Oun\sum_{n}e^{2\varepsilon n}m\big(R\ge n\big)
$$
and this quantity is finite if $\varepsilon$ is small enough by
using (\ref{EXPRETURN}) and property (I)-(ii) in \cite[Section
3.2]{young}.
Collecting all the bounds we get the following upper-bound for the
first term in (\ref{finfin}):
$$
2\sum_{k=1}^{n-1}\ \int\ \Bigl ( {\pf}^k\
U(y,\ldots,{\fbar}^{n-k-1}(y))-\ {\pf}^{k+1}\ U(
{\fbar} (y),\ldots,{\fbar}^{n-k-1}(y))\Bigr )^2\ d\mubar(y) \leq
$$ 
$$
\Oun \ \sum_{k=1}^{\infty}\left(\sum_{l=1}^{k-1}\Gamma^{k}_{l}\right)^{2}\,.
$$
Choosing $p_0$ large enough in the definition (\ref{lulu}) of $\ell$, we have 
$$
B_l \leq L_{l} (\textup{diam}(M))^\eta + C^{\eta} \sum_{j=l}^{\infty} \frac{L_{j+1}}{(j-l+1)^2}\,\cdot
$$
Using several times Young's inequality, one easily gets
$$
\sum_{k=1}^{\infty}\left(\sum_{l=1}^{k-1}\Gamma^{k}_{l}\right)^{2}\le \Oun \sum_{j=1}^{n}L_{j}^{2}\;.
$$
Since a separately H\"older continuous function of $n$ variables can also be
considered as a separately H\"older continuous function of $n+k$
($k>0$) with $L_j=0$ for $j>n$, the same estimate holds for the second
term in (\ref{finfin}).

We now prove that the third term in (\ref{finfin}) tends to zero when
$p\to\infty$. 
Using Lemma \ref{wreg} and estimate (\ref{salade}) with $k=n$ and $l=n-1$
we observe that $S_n=\pf^n U_n$ belongs to the Banach space $X_\sigma$.
The result follows at once using Lemma \ref{LYI}.

This ends the proof of Theorem \ref{Main}.


\end{document}